\documentclass[reqno]{amsart}
\usepackage{amssymb,graphics}

\newcommand{\FF} {{\mathbb F}}

\newcommand{\CCC}{{\mathcal C}}
\newcommand{\QQ} {{\mathbb Q}}
\newcommand{\NN} {{\mathbb N}}
\newcommand{\ZZ} {{\mathbb Z}}
\newcommand{\im} {{\mathop{\rm im}\nolimits}}
\newcommand{\rank} {{\mathop{\rm rank}\nolimits}}
\newcommand{\simplicial} {{\mathcal S}}
\newcommand{\comment}[1]{}

\newtheorem{theorem}{Theorem}
\newtheorem{corollary}[theorem]{Corollary}
\newtheorem{proposition}[theorem]{Proposition}

\newtheorem{lemma}[theorem]{Lemma}
\theoremstyle{remark}
\newtheorem{remark}[theorem]{Remark}

\begin{document}
\title{Cyclotomic and simplicial matroids}

\author{Jeremy L. Martin}
\address{School of Mathematics\\
University of Minnesota\\
Minneapolis, MN 55455}
\email{martin@math.umn.edu}

\author{Victor Reiner}
\address{School of Mathematics\\
University of Minnesota\\
Minneapolis, MN 55455}
\email{reiner@math.umn.edu}


\thanks{First author supported by NSF Postdoctoral Fellowship.
Second author supported by NSF grant DMS-0245379.}

\begin{abstract}
Two naturally occurring matroids representable over $\QQ$ are
shown to be dual:  the {\it cyclotomic matroid} $\mu_n$ represented
by the $n^{th}$ roots of unity $1,\zeta,\zeta^2,\ldots,\zeta^{n-1}$
inside the cyclotomic extension
$\QQ(\zeta)$,
and a direct sum of copies of a certain simplicial
matroid, considered originally by Bolker in the context
of transportation polytopes.
A result of Adin leads to an upper bound for the number of 
$\QQ$-bases for $\QQ(\zeta)$ among the $n^{th}$ roots of unity, which
is tight if and only if $n$ has at most two odd prime factors.
In addition, we study the Tutte polynomial of $\mu_n$ in the case that $n$
has two prime factors.
\end{abstract}

\maketitle

\section{Introduction}

This paper is about two matroids representable over $\QQ$ that
turn out, somewhat unexpectedly, to be {\it dual} (or {\it orthogonal}).
Briefly, a matroid is a combinatorial abstraction of the linear
dependence data associated
to a (finite) set of vectors in a vector space.  That is,
the data for a matroid on ground set $E$ records which subsets of $E$
are dependent and independent, or the linear span and dimension of each subset, etc.
The matroid is {\it representable} over a field $\FF$ if the elements of
$E$ may be identified
with vectors in an $\FF$-vector space that achieve the matroid data.
The {\it dual} of a matroid $M$ on $E$ is defined generally as
the matroid $M^*$ whose bases are the complements of bases in $M$; this abstracts
the situation where $M, M^*$ are matroids represented over $\FF$ by the columns of
two matrices whose row spaces are orthogonal complements of each other.
For background, motivation, and matroid terminology left undefined here, 
see any of the standard references 
\cite{CrapoRota, Oxley, Welsh, White1, White2, White3}.

The first matroid that we study relates to cyclotomic extensions.  
Let $\zeta_n$ be a primitive $n^{th}$ root of unity, which we will
abbreviate by $\zeta$ when no confusion can arise.
Recall that the degree of the cyclotomic extension $\QQ(\zeta)$ 
is given by Euler's $\phi$-function
$\phi(n)$, with the following formula:  
if $n=p_1^{m_1} \cdots p_r^{m_r}$ for distinct primes $p_i$, then
$$
\phi(n) = p_1^{m_1-1} \cdots p_r^{m_r-1} (p_1-1) \cdots (p_r-1).
$$
The $\QQ$-vector space $\QQ(\zeta)$ gives rise to a natural matroid
$\mu_n$, represented
by the $n$ vectors $Z_n=\{1,\zeta,\zeta^2, \ldots, \zeta^{n-1}\}$ and
having rank $\phi(n)$.  We call $\mu_n$ the
{\it cyclotomic matroid of order~$n$}.

The second matroid is a {\it simplicial matroid}~\cite[Chapter 6]{White2};
that is, it is represented by the columns of a simplicial boundary
map (see, e.g.,~\cite{Munkres} for background on simplicial homology).
Let $\Delta=\Delta^d$ be a $d$-dimensional simplicial
complex\footnote{We adopt the following notational convention throughout.
When we wish to emphasize that a simplicial complex $\Delta$ has dimension
$d$, we will denote it by the symbol $\Delta^d$; otherwise, we will
frequently omit the superscript to simplify notation.}
whose (reduced) simplicial homology with coefficients in $\FF$ 
vanishes in codimension~$1$; that is, 
$\tilde{H}_{d-1}(\Delta,\FF)=0$.  If $\Delta$ has exactly $n$ facets (that
is, faces of maximum dimension~$d$),
then the $n$ columns in the matrix expressing the simplicial boundary
map 
$$
\tilde{C}_d(\Delta^d,\FF) \overset{\partial_d}{\longrightarrow} \tilde{C}_{d-1}(\Delta^d,\FF)
$$ 
represent a matroid over $\FF$ that we will call the {\it simplicial 
matroid} $\simplicial(\Delta^d,\FF)$.  
This is a matroid on $n$ elements, with
$$
\rank\;\simplicial(\Delta^d,\FF) ~=~
\dim_\FF\,\im\:\partial_d ~=~ \dim_\FF\,\ker\,\partial_{d-1}.
$$

In this paper, we consider the $(r-1)$-dimensional 
simplicial complex $\Delta^{r-1}_{n_1,\ldots,n_r}$ defined as
the {\it simplicial join} \cite[\S62]{Munkres}
of $0$-dimensional complexes $\Delta^0_{n_1},
\ldots, \Delta^0_{n_r}$,
where $\Delta^0_{n_i}$ consists of $n_i$ disjoint
vertices.  That is, a face of $\Delta^{r-1}_{n_1,\ldots,n_r}$
contains at most one vertex from each $\Delta^0_{n_i}$.
Our main result expresses the connection between cyclotomic matroids
and these particular simplicial matroids.

\begin{theorem}
\label{main-theorem}
Let $n = p_1^{m_1} \cdots p_r^{m_r}$, with $p_1,\ldots,p_r$  distinct primes 
and $m_1,\ldots,m_r$ positive integers.

Then the following two matroids representable over $\QQ$ are dual:
\begin{itemize}
\item The cyclotomic matroid $\mu_n$.
\item The direct sum of
$p_1^{m_1-1} \cdots p_r^{m_r-1}$ copies of
$\simplicial(\Delta^{r-1}_{p_1,\ldots,p_r},\QQ)$.  
\end{itemize}
\end{theorem}


\noindent
Bolker~\cite{Bolker} was the first to study the simplicial matroid 
$\simplicial(\Delta^{r-1}_{n_1,\ldots,n_r},\QQ)$, where the $n_i$ are
positive integers, not necessarily prime.  He proposed the
bases of this matroid as higher-dimensional analogues of
spanning trees\footnote{This is not the only way to generalize
the notion of spanning tree to 2- or higher-dimensional simplicial complexes;
see, e.g.,~\cite{BeinekePippert}, \cite{Faridi}, \cite{HararyPalmer}.  These other
generalizations, however, play no role in our present study.}
in bipartite graphs,
and studied their relation to vertices of certain transportation
polytopes.  Each such basis $T$ of $\simplicial(\Delta^{r-1}_{n_1,\ldots,n_r},\QQ)$
gives rise to a $\QQ$-acyclic pure $(r-1)$-dimensional
simplicial complex $\Delta_T$, obtained by attaching the $(r-1)$-simplices
indexed by $T$ to the $(r-2)$-skeleton of $\Delta^{r-1}_{n_1,\ldots,n_r}$.
He also suggested
\begin{equation}
\label{Bolker-bound}
\prod_{i=1}^r n_i^{\prod_{j \neq i}(n_j-1)}
\end{equation}
as an upper bound for the number of bases,
and showed that this formula gives
the exact number of bases if and only if
at most two of the $n_i$ exceed $2$.  In particular, Bolker
showed~\cite[Theorems~27,~28]{Bolker} that
this occurs exactly when for each basis $T$, the (pure torsion) integer
homology group $\tilde{H}_{d-1}(\Delta_T,\ZZ)$ is trivial.

A result of Adin~\cite{Adin} (generalizing work of Kalai in~\cite{Kalai})
implies that \eqref{Bolker-bound} is indeed an 
upper bound for the number of bases of $\Delta^{r-1}_{n_1,\ldots,n_r}$.
The method of Adin and Kalai is to generalize the
classical Kirchhoff Matrix-Tree Theorem to simplicial complexes, using
the Binet-Cauchy determinant formula: this yields the result
\begin{equation}
\label{Adin-formula}
    \sum_{\substack{\text{bases } T \text{ of }\\ 
           \simplicial( \Delta^{r-1}_{n_1,\ldots,n_r},\QQ)} }
    \left\vert \tilde{H}_{r-2}(\Delta_T,\ZZ) \right\vert^2 ~=~
    \prod_{i=1}^r n_i^{\prod_{j \neq i}(n_j-1)}.
\end{equation}
That is, Bolker's formula~\eqref{Bolker-bound} is an {\it exact} count
for such bases with the following weighting:
instead of a basis $T$ contributing $1$ to the count,
its contribution is the square of the order
of its integral homology group in codimension~$1$.

Combining the results of Adin and Bolker with Theorem~\ref{main-theorem}, 
and the fact that dual matroids have the same number of bases,
one immediately obtains the following corollary.

\begin{corollary}
\label{upper-bound-corollary}
Let $n= p_1^{m_1} \cdots p_r^{m_r}$, with $p_1,\ldots,p_r$  distinct primes 
and $m_1,\ldots,m_r$ positive integers, 
and let $\zeta$ be a primitive $n^{th}$ root of unity.

Then the number of subsets of $Z_n=\{1,\zeta,\zeta^2,\dots,\zeta^{n-1}\}$
that are bases for $\QQ(\zeta)$ is bounded above by
$$
\left(\:
   \prod_{i=1}^r p_i^{\prod_{j \neq i}(p_j-1)} 
\right)^{p_1^{m_1-1} \cdots p_r^{m_r-1}},
$$
with equality if and only if $n=2^a p^b q^c$ for
nonnegative integers $a,b,c$.
\end{corollary}

Section~\ref{proof-section} contains the proof of Theorem~\ref{main-theorem}.
In Section~\ref{two-primes-section},
we study the case that $n=p_1^{m_1}p_2^{m_2}$
has only two prime factors.  Here the
simplicial complex $\Delta^1_{p_1,p_2}$ is just the
complete bipartite graph $K_{p_1,p_2}$, and the matroid
$\simplicial(\Delta^1_{p_1,p_2},\QQ)$ is the {\it cycle matroid} (or {\it graphic
matroid}) of $K_{p_1,p_2}$, whose bases are spanning trees.
We study enumerative invariants finer than the number of bases for
these particular graphic matroids,
such as their Tutte polynomials, and a 
weighted enumeration of their spanning forests.


\section{Proof of Theorem~\ref{main-theorem}}
\label{proof-section}

We first state a well-known general fact about matroids, to be used in the proof.

\begin{lemma}
\label{general-matroid-lemma}
Let $M$ be a matroid of rank $\rho$ on ground set $E$.
Suppose there exists a (disjoint) decomposition $E = \bigsqcup_{i=1}^t E_i$
such that $\rank(M) = \sum_{i=1}^t \rank\;(M|_{E_i})$.

Then $M \cong \bigoplus_{i=1}^t M|_{E_i}$.
\end{lemma}

\begin{proof}
Because of the assumption on ranks, a basis for $M$
must be a disjoint union of bases of the matroids $M|_{E_i}$.
\end{proof}

\vskip .2in
\noindent {\it Proof of Theorem~\ref{main-theorem}.}  

%
We first reduce to the case where $n$ is square-free.
Let $s=p_1 \cdots p_r$ denote the square-free part of $n$,
and let
$$
t=\frac{n}{s} =p_1^{m_1-1} \cdots p_r^{m_r-1}.
$$
Note that 
$$
\rank\;\mu_n = \phi(n) = t \cdot \phi(s) = t \cdot \rank\;\mu_s.
$$
For $0 \leq j \leq t-1$, let $E_j :=
\{\zeta^j, \zeta^{j+t}, \dots, \zeta^{j+(s-1)t}\}$.
Then
$$
Z_n = \bigsqcup_{j=0}^{t-1} E_j,
$$
and the submatroid $\mu_n|_{E_j}$ obtained by restricting
$\mu_n$ to the ground set $E_j$ is isomorphic
to $\mu_s$.
By Lemma~\ref{general-matroid-lemma}, it follows that
$\mu_n$ is isomorphic to a direct sum of $t$ copies of $\mu_s$.
The following figure illustrates the decomposition for $n=18$,
$s=6$ and $t=3$.

\vskip 0.2in
\begin{center}
\resizebox{4.5in}{1.5in}{\includegraphics{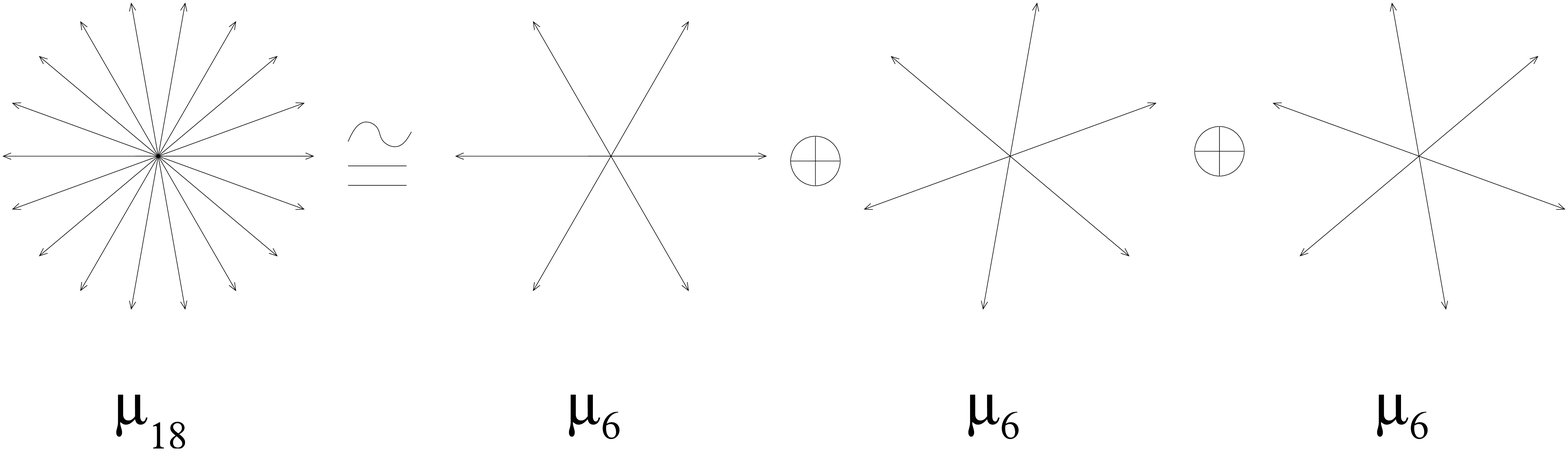}}
\end{center}
\vskip 0.2in

\noindent
Since duality commutes with direct sums, 
we have reduced to the case that
$n=p_1 \cdots 
p_r$
is square-free, and assume this for the remainder of the proof.

For each prime $p$ dividing $n$, consider
the short exact sequence of $\QQ$-vector spaces
$$
0 \rightarrow \QQ \overset{i_p}{\longrightarrow} 
              \QQ^p \overset{\pi_p}{\longrightarrow} \QQ(\zeta_p) \rightarrow 0
$$
in which $\pi_p$ maps the $j^{th}$ standard basis element $e_j$ in
$\QQ^p$ to $\zeta_p^j$, and the image of $i_p$ is the 
line $\QQ(e_1 + \cdots + e_n)$.  Rephrasing this, the cochain 
complex $\CCC(p)$ defined by
$$
\begin{matrix}
0 & \rightarrow & \CCC^0(p) & \longrightarrow               & \CCC^1(p)
  & \rightarrow &0 \\
0 & \rightarrow & \QQ       & \overset{i_p}{\longrightarrow} & \QQ^p     
  & \rightarrow &0 \\
\end{matrix}
$$
has 
$$
\begin{aligned}
H^0(\CCC(p)) &=0 \\
H^1(\CCC(p)) &=\QQ(\zeta_p) \,\, (\cong \QQ^{p-1}).
\end{aligned}
$$
Note that the complex $\CCC(p)$ coincides, up to a shift in homological degree,
with the augmented simplicial cochain complex that computes
the cohomology of the $0$-dimensional complex $\Delta^0_p$.
Similarly, the tensor product of complexes 
$$
\CCC:=\bigotimes_{j=1}^r \CCC(p_j)
$$
coincides with the augmented simplicial cochain complex for the
simplicial join
$\Delta_{p_1,\ldots,p_r} \cong \Delta^0_{p_1} * \cdots * \Delta^0_{p_r}$.
The K\"unneth formula
over $\QQ$ then implies that the complex $\CCC$
has no homology except in the last position.  The last cochain group
in $\CCC$ maps surjectively onto the homology as follows:
\begin{equation}
\label{tensor-cokernel}
\CCC^r \ := \ \bigotimes_{j=1}^r \CCC^1(p_j) \ \cong \
\bigotimes_{j=1}^r \QQ^{p_j} \ 
\xrightarrow{\otimes_j \pi_{p_j}} \
\bigotimes_{j=1}^r \QQ(\zeta_{p_j}). 
\end{equation}
The Chinese Remainder Theorem  allows one to
identify the last map in \eqref{tensor-cokernel}
with the surjection $\QQ^n \overset{\pi_n}{\rightarrow} \QQ(\zeta)$
sending $e_j$ to $\zeta^j$, that is, the map whose matrix has
columns representing $\mu_n$.  Because $\CCC$ has no other homology,
the last coboundary map
$\delta: \CCC^{r-1} \rightarrow \CCC^{r}$ has columns spanning the
kernel of $\pi_n$.  Its transpose 
$$
C_{r-1}(\Delta_{p_1,\ldots,p_r}; \QQ) 
\ \overset{\delta^T}{\longrightarrow} \ 
C_{r-2}(\Delta_{p_1,\ldots,p_r}; \QQ)
$$ 
has columns which represent the simplicial
matroid $\simplicial(\Delta_{p_1,\ldots,p_r}, \QQ)$, and hence this
matroid is dual to $\mu_n$.
\qed
\vskip .2in

We close this section with a few remarks.


\vskip .2in
\begin{remark}
To give some further context for Theorem~\ref{main-theorem}, it should
be noted that matroid duality has played a key role
in the study of simplicial matroids from the beginning.
Crapo and Rota studied the full simplicial matroid 
$\simplicial^m_{m-k}(\FF)$ corresponding to the $(k-1)$-skeleton of 
the full simplex with $m$ vertices, and observed
that Alexander duality implies that $\simplicial^m_k(\FF)$ and
$\simplicial^m_{m-k}(\FF)$ are dual matroids \cite[Proposition 11.4]{CrapoRota} 
(see also \cite[Theorem 6.2.1]{White2}).  Lindstr\"om \cite{Lindstrom} observed another
interesting example of such a 
duality, this time involving two natural matroids
(one of them simplicial), both of
whose ground sets are indexed by the bases of some
(other) simple matroid; see~\cite[\S 6.5]{White2}.
\end{remark}

\vskip .2in
\begin{remark}
It has been observed (e.g., by Johnsen \cite{Johnsen}) 
that the primitive $n^{th}$ roots of unity form a $\QQ$-basis for
$\QQ(\zeta)$ if and only if $n$ is square-free.  Thus if
$n=p_1 \cdots p_r$
is square-free, one might ask which particular basis of the simplicial matroid
$\simplicial(\Delta^{r-1}_{p_1,\ldots,p_r},\QQ)$ (that is, which
higher-dimensional tree $T$ of $\Delta^{r-1}_{p_1,\ldots,p_r}$) 
corresponds to the set of primitive $n^{th}$ roots.
Passing through matroid duality, and the Chinese Remainder Theorem
identification in the proof of the theorem, 
one can check that $T$ is the union of the stars within
$\Delta^{r-1}_{p_1,\ldots,p_r}$ of the vertices
$v_1,\dots,v_r$, where $v_i$ is the unique vertex of $\Delta^0_{p_i}$ 
labeled by $0$ modulo $p_i$.
Applying a simple nerve argument to this cover of $T$ by stars,
one finds that $T$ is not only $\QQ$-acyclic as expected, but in fact
contractible.  Hence it has vanishing homology group $\tilde{H}_{r-2}(\Delta_T,\ZZ)$,
and contributes exactly $1$ in the summation of equation \eqref{Adin-formula}.
\end{remark}

\vskip .2in
\begin{remark}
If $n$ is divisible by at most two primes (that is, $r\leq 2$ above),
then the complex $\Delta^{r-1}_{p_1,p_2}$ is a graph.
In this case, the matroid $\mu_n$ is {\it cographic} (see,
e.g.,~\cite[p.~35]{White2}, and Section~\ref{two-primes-section} below).  

Furthermore, if $n$ is odd,
then the cyclotomic
matroid $\mu_{2n}$ is the parallel extension of $\mu_n$ in which one
creates one parallel copy for each ground set element.  
For instance, $\mu_n$ is again cographic for $n=2pq$, where $p,q$ are
odd primes: it is the cographic matroid of the graph obtained by duplicating
every edge of the complete bipartite graph $K_{p,q}$.
\end{remark}

\section{The case of two prime factors}
\label{two-primes-section}
The goal of this section is to compute some enumerative invariants of the
cyclotomic matroid $\mu_n$, finer than the number of bases,
in the special case where $n=p_1^{m_1}p_2^{m_2}$
has only two prime factors.

The {\it Tutte polynomial} $T_M(x,y)$ is one of the most important
isomorphism
invariants of a matroid $M$; see the excellent survey by Brylawski and
Oxley~\cite[Chapter 6]{White3} for background.  We begin by
reviewing briefly the definition of the Tutte polynomial,
as well as some elementary properties relevant for our calculations.

Let $M$ be a matroid with ground set $E$, and let $r$ be
the rank function on subsets $A$ of $E$.  When $M$ is representable,
$r(A)$ is the dimension of the linear span of the vectors in $A$.
The {\it Tutte polynomial} $T_M(x,y)$ may be defined as
the corank-nullity generating function
    \begin{equation} \label{Tutte-definition}
      T_M(x,y) := \sum_{A \subset E} (x-1)^{r(M)-r_M(A)}\,
        (y-1)^{|A|-r_M(A)}
    \end{equation}
\cite[\S6.2]{White3}.  Two easy consequences
of~\eqref{Tutte-definition} are as follows.
First, if $M^*$ is the dual matroid to $M$, then
  \begin{equation} \label{dual-tutte}
  T_{M^*}(x,y)=T_M(y,x).
  \end{equation}
Second, setting $y=1$ and replacing $x$ with $x+1$
in \eqref{Tutte-definition} gives a
generating function for independent subsets of $M$ according to
their cardinality:
  \begin{equation} \label{indep-set-gf}
  T_M(x+1,1)=\sum_{\text{independent }I \subset E} x^{r(M)-|I|}.
  \end{equation}

Let $G=(V,E)$ be a graph, and let $M=M(G)$ be the
corresponding graphic matroid on ground set $E$,
whose bases are the spanning trees of $G$ 
and whose independent sets are acyclic subgraphs (= forests).
Then there is a substitution of variables in 
$T_G:=T_{M(G)}$ which gives Crapo's {\it coboundary polynomial}
$\bar\chi_G(q,t)$ \cite[\S 6.3F]{White3}:
  \begin{equation}
  \label{Crapo-defn}
  \begin{aligned}
  \bar\chi_G(q,t) &~:=~
    q^{-1}\sum_{\substack{\text{ vertex colorings }
                    \\f:V \rightarrow \{1,2,\ldots,q\}}} 
     t^{\left|\{\{v_1,v_2\} \in E:\; f(v_1)=f(v_2)\}\right|} \\
     & ~=~ (t-1)^{|V|-c(G)}\: T_G\left(\frac{q+t-1}{t-1},\;t\right),
  \end{aligned}
  \end{equation}
where $c(G)$ denotes the number of connected components of $G$.
Note that $T_G(x,y)$ may be recovered from $\bar\chi_G(x,y)$
via the substitution $q=(x-1)(y-1),t=y$, and that the usual
chromatic polynomial $\chi_G(q)$ is the specialization $q \, \bar\chi_G(q,0)$.

In this section, we consider the special case that
$n=p_1^{m_1} p_2^{m_2}$.  
By Theorem~\ref{main-theorem}, the cyclotomic matroid
$\mu_n$ is dual
to the direct sum of $p_1^{m_1-1} p_2^{m_2-1}$ copies of the 
simplicial matroid $\simplicial(\Delta^1_{p_1,p_2};\QQ)$.
As mentioned previously,
the $1$-dimensional complex $\Delta^1_{p_1,p_2}$ is simply
the complete bipartite graph $K_{p_1,p_2}$, and
the simplicial matroid is the usual graphic matroid
$M(K_{p_1,p_2})$.

We begin by obtaining an exponential generating function
for the coboundary polynomials of these cycle matroids,
mimicking Ardila's method in~\cite[Theorem 2.4.1]{Ardila}.
Let $K_{p_1,p_2}$ have bipartite
vertex set $V_1 \sqcup V_2$ with $|V_i|=p_i$.  Every vertex-coloring 
$f:V \rightarrow \{1,2,\ldots,q\}$ decomposes the
partite sets $V_1,V_2$
into (possibly empty) color classes $V_1^{(i)}$, $V_2^{(i)}$:
        $$V_1 = \bigsqcup_{i=1}^q V_1^{(i)}, \qquad\qquad
        V_2 = \bigsqcup_{i=1}^q V_2^{(i)}$$
with
$$
|\{\{v_1,v_2\} \in E:\; f(v_1)=f(v_2)\}| 
 ~=~ \sum_{i=1}^q |V_1^{(i)}| |V_2^{(i)}|.
$$
Now exponential generating function manipulation
(see~\cite[Prop.~5.1.3]{Stanley}) gives the following formula.

\begin{proposition} \label{Kab-EGF}
Denote by $\bar{\chi}_{p_1,p_2}(q,t)$
the coboundary polynomial of the graphic matroid of the complete
bipartite graph $K_{p_1,p_2}$.  Then
$$
 1+q \left(
      \sum_{\substack{(p_1,p_2) \in \\ \NN^2-\{(0,0)\}}}
         \bar{\chi}_{p_1,p_2}(q,t) \;
          \frac{x_1^{p_1} x_2^{p_2}}{p_1!\,p_2!}
        \right) \quad =  \quad
      \left( \sum_{(m_1,m_2) \in \NN^2} t^{m_1 m_2} 
        \frac{x_1^{m_1} x_2^{m_2}}{m_1!\,m_2!} \right)^q.
$$
\end{proposition}

\noindent
This formula does not appear to
generalize to the case that $n$ has more than two prime factors.
\comment{  The correct
generalization of the coboundary polynomial seems to be
        $$\bar{\chi}_{p_1,\dots,p_r}(q,t) ~=~
    q^{n-\phi(n)-|V|}\sum_{\substack{\text{ vertex colorings }
                    \\f:V \rightarrow \{1,2,\ldots,q\}}} 
     t^{\left|\{\text{facets }\{v_1,\dots,v_r\}:\; \sum_i (-1)^if(v_i)=0\}\right|}$$
but it is unclear how to obtain an exponential generating function for these
polynomials.
}

\begin{remark}
We remark that setting $t=0$ yields
the following simple exponential generating function for
the chromatic polynomial $\chi_{p_1,p_2}(q)$ of $K_{p_1,p_2}$:
$$
 1+ \sum_{(p_1,p_2) \in  \NN^2-\{(0,0)\}}
         \,\chi_{p_1,p_2}(q) \;
          \frac{x_1^{p_1} x_2^{p_2}}{p_1!\,p_2!}
         \quad =  \quad
      \left(  e^{x_1} + e^{x_2} - 1 \right)^q.
$$
\end{remark}

We now apply the generating function for coboundary polynomials
to enumerate the $\QQ$-linearly independent subsets
of the $n^{th}$ roots of unity.
Setting $t=y+1$ and $x_i=z_i/y$ in
the formula of Proposition~\ref{Kab-EGF}
and applying~\eqref{Crapo-defn} gives
    \begin{equation} \label{Tutte-EGF-before-limit}
     \begin{aligned}
      &y^{-1} \ 
        \sum_{(p_1,p_2) \in \NN^2-\{(0,0)\}}
           T_{K_{p_1,p_2}}\left(\frac{q+y}{y},\:y+1 \right)
            \frac{z_1^{p_1} z_2^{p_2}}{p_1!\,p_2!} \quad = \\
     & \qquad \qquad
      q^{-1} \left[ \left( \sum_{(m_1,m_2) \in \NN^2} 
             \frac{(y+1)^{m_1 m_2}}{y^{m_1+m_2}}\:
          \frac{z_1^{m_1} z_2^{m_2}}{m_1!\,m_2!} \right)^q - 1 \right].
     \end{aligned}
    \end{equation}

\noindent
Multiplying~\eqref{Tutte-EGF-before-limit} through by $y$
and taking the limit
as $q$ approaches $0$ (via L'H\^opital's Rule), we obtain
  \begin{equation} \label{cyclotomic-indep-EGF}
    \begin{aligned}
        &\sum_{(p_1,p_2) \in \NN^2-\{(0,0)\}}
           T_{K_{p_1,p_2}}(1,y+1)~
            \frac{z_1^{p_1} z_2^{p_2}}{p_1!\,p_2!} \quad = \\
      & \qquad \qquad 
        z \log \left( \sum_{(m_1,m_2) \in \NN^2}
           \frac{(y+1)^{m_1 m_2}}{y^{m_1+m_2}}~
          \frac{z_1^{m_1} z_2^{m_2}}{m_1!\,m_2!} \right).
     \end{aligned}
    \end{equation}

\noindent
Once again, recall that for $n=p_1^{m_1}p_2^{m_2}$,
Theorem~\ref{main-theorem} implies that the cyclotomic matroid
$\mu_n$ is dual
to the direct sum of $p_1^{m_1-1} p_2^{m_2-1}$ copies of
the graphic matroid  $M(K_{p_1,p_2})$.  Therefore,~\eqref{dual-tutte}
and~\eqref{indep-set-gf} imply the following result.

\vskip .1in
\begin{corollary}
\label{indep-set-gf-corollary}
Let $n=p_1^{m_1} p_2^{m_2}$, where $p_1,p_2$ are distinct primes.

Then the generating function
    $$
    \sum_{\substack{I \subset Z_n \\ \QQ{\rm-linearly~independent}  }} y^{\phi(n)-|I|}
    $$
equals the $(p_1^{m_1-1} p_2^{m_2-1})^{th}$ power
of the coefficient of $\frac{z_1^{p_1} z_2^{p_2}}{p_1!\,p_2!}$ in 
    $$
    y \log \left( \sum_{(m_1,m_2) \in \NN^2} 
        \frac{(y+1)^{m_1 m_2}}{y^{m_1+m_2}}~
        \frac{z_1^{m_1} z_2^{m_2}}{m_1!\,m_2!} \right).
    $$
\end{corollary}
\vskip .1in

We conclude this section with an observation about the independent set
polynomial $T_M(x+1,1)$ (see~\eqref{indep-set-gf}) for $M=M(K_{p,q})$.
Note that $T_M(x+1,1)$ 
enumerates spanning subsets of $Z_{pq}$ by cardinality 
when $p,q$ are distinct primes.
Computations
for small values of $p,q$ suggest that for $p\leq q$ (not necessarily prime), 
$T_{M(K_{p,q})}(x+1,1)$ is divisible by
$(x+p)^{q-p+1}$.  This follows from a more refined statement
(Proposition~\ref{more-refined-statement} below) 
on the enumeration of spanning forests in $K_{p,q}$.
Let $K_{p,q}$ have disjoint vertex sets 
$V=\{v_1,\dots,v_p\}$ and $W=\{w_1,\dots,w_q\}$.
For a subset $F$ of its edges and $u$ a vertex,
let $\deg_F(u)$ denote the degree of $u$ in the edge-subgraph corresponding
to $F$.  Define a forest enumerator
$$
A_{p,q} ~=~ A_{p,q}(x_1,\dots,x_p) ~:=~ \sum_{\text{forests}\;F}
        x_1^{\deg_F(v_1)} \dots x_p^{\deg_F(v_p)}.
$$

For example, it can be calculated that 
$$
\begin{aligned}
A_{1,q} &= (1+x_1)^q, \\
A_{2,q} &= (1+x_1+x_2)^{q-1}(1 + x_1 + x_2 + q x_1 x_2).
\end{aligned}
$$

\noindent
Note that $A_{p,q}$ specializes to the independent set
polynomial of $K_{p,q}$ as follows.  Since $K_{p,q}$ has $pq$ edges,
each of which is incident to exactly one $v_i$, we obtain
        $$T_{K_{p,q}}(x+1,1)\;=\;x^{pq}A_{p,q}(x^{-1},\,\dots,\,x^{-1}).$$

\vskip .1in
\begin{proposition}
\label{more-refined-statement}
Let $q\geq p \geq 1$ be integers.  Then 
$$
A_{p,q} = \sum_{j=0}^{p-1} (1+x_1+\dots+x_p)^{q-j} \tilde{A}_{p,j}
$$
where $\tilde{A}_{p,q}:=\sum_{F} x_1^{\deg_F(v_1)} \dots x_p^{\deg_F(v_p)}$,
in which the summation is over spanning forests $F$ in $K_{p,q}$
in which every vertex $w_i$ has degree at least $2$.

In particular, $A_{p,q}$ is divisible by $$(1+x_1+\dots+x_p)^{q-p+1}.$$
\end{proposition}
\vskip .1in

\noindent

\vskip .1in
\begin{proof}
Given  $F$ a spanning forest of $K_{p,q}$, let $W' \subset W$
be the set of $w_i$ for which $\deg_F(w_j)>1$.  We claim that $|W'| <p $,
else $F$ contains at least $2p$ edges among the $2p$ vertices
$|W' \cup V|$, which contradicts the assumption that $F$ is acyclic.

After choosing the restriction of $F$ to the vertices $W' \cup V$,
the remaining edges of $F$ consist of at most one edge incident to
each vertex in $W-W'$.  Classifying $F$ according to the cardinality
$|W'|=j$, the assertion follows.
\end{proof}

The quotient $\frac{A_{p,q}}{(1+x_1+\dots+x_p)^{q-p+1}}$ 
is a symmetric polynomial of total degree $2(p-1)$
in the variables $x_1,\dots,x_p$.  It appears not to factor further in general.
One might also hope for a factorization of a forest enumerator that keeps track of
the degrees in $V$ as well as those in $W$, but this also appears not to factor further.


\begin{thebibliography}{99}

\bibitem{Ardila}
F. Ardila,
Enumerative and algebraic aspects of matroids and hyperplane arrangements.
Ph.D. thesis, MIT, 2003.

\bibitem{Adin}
R.M. Adin, 
Counting colorful multi-dimensional trees. 
\textit{Combinatorica} {\bf 12} (1992), 247--260.

\bibitem{BeinekePippert}
L.W. Beineke and R.E. Pippert,
Properties and characterizations of $k$-trees.
\textit{Mathematika} {\bf 18} (1971), 141--151.

\bibitem{Bolker}
E.D. Bolker, 
Simplicial geometry and transportation polytopes. 
\textit{Trans. Amer. Math. Soc.} {\bf 217} (1976), 121--142.

\bibitem{CrapoRota}
H.H. Crapo and G.-C. Rota,
On the foundations of combinatorial theory: Combinatorial geometries. 
Preliminary edition. The M.I.T. Press, Cambridge, Mass.-London, 1970.

\bibitem{Faridi}
S. Faridi,
The facet ideal of a simplicial complex.
\textit{Manuscripta Math.} {\bf 109} (2002), no. 2, 159--174.

\bibitem{HararyPalmer}
F. Harary and E.M. Palmer,
On acyclic simplicial complexes.
\textit{Mathematika} {\bf 15} (1968), 115--122.

\bibitem{Johnsen}
K. Johnsen,
Lineare Abh\"angigkeiten von Einheitswurzeln.
\textit{Elem. Math.} {\bf 40} (1985), 57--59.

\bibitem{Kalai}
G. Kalai, 
Enumeration of $\QQ$-acyclic simplicial complexes. 
\textit{Israel J. Math.} {\bf 45} (1983), 337--351.

\bibitem{Lindstrom}
B. Lindstr\"om,
Matroids on the bases of simple matroids.
\textit{Europ. J. Combin.} {\bf 2} (1981), 61--63.

\bibitem{Munkres}
J.R. Munkres,
Elements of algebraic topology. 
Addison-Wesley Publishing Company, Menlo Park, CA, 1984.

\bibitem{Oxley}
J.G. Oxley,
Matroid theory.
Oxford Science Publications. 
The Clarendon Press, Oxford University Press, New York, 1992.

\bibitem{Stanley}
R.P. Stanley,
Enumerative Combinatorics, vol. 2.
Cambridge University Press, Cambridge, 1999.

\bibitem{Welsh}
D.J.A. Welsh,
Matroid theory. 
London Math. Soc. Monographs {\bf 8}
Academic Press, London-New York, 1976.

\bibitem{White1}
N. White,
Theory of matroids. 
\textit{Encyclopedia of Mathematics and its Applications} {\bf 26}. 
Cambridge University Press, Cambridge, 1986.

\bibitem{White2}
N. White,
Combinatorial geometries. 
\textit{Encyclopedia of Mathematics and its Applications} {\bf 29}. 
Cambridge University Press, Cambridge, 1987.

\bibitem{White3}
N. White,
Matroid applications.
\textit{Encyclopedia of Mathematics and its Applications} {\bf 40}.
Cambridge University Press, Cambridge, 1992. 

\end{thebibliography}
\end{document}